\mag=\magstep 1
\baselineskip=15pt plus .1 pt 
\hsize=17.0truecm
\vsize=25.0truecm

\footline{\ifnum\pageno=1{ \hss\tenrm-- 1 --\hss } \else \arith \fi}
            \def \arith{ \hss\tenrm-- \folio\ --\hss }
\pageno=1

\voffset=-.5truecm
\hoffset=-.5truecm

\def\square{\vcenter{\vbox {\hrule \hbox{\vrule height 6pt
                                    \kern 6pt\vrule }\hrule }}}
\def\tmu{T_{\cal M}^{\theta } }
\def\i=1{\, i=1,\dots ,n\, }

\def\nore{\| \sum _{i=1}^m e_i \| }
\def\ncl{ {1\over {n^{{1\over q}} } } }

\centerline {{\bf BANACH SPACES OF THE TYPE OF TSIRELSON}}
\centerline {{\bf  S. A. Argyros and I. Deliyanni}}
\medskip
 {\bf Abstract }. {\sl { To any pair $\, (\cal M \, ,\, \theta
\,)\, $, where $\, \cal M \,$ is a family of finite subsets of $\, {\bf N }\, ,$
compact in the pointwise topology, and $\, 0<\theta < 1 \, ,$ we associate a
Tsirelson-type Banach space $\, T_{\cal M}^{\theta }\, .$ It is shown that if
the Cantor-Bendixson index of $\, \cal M \, $ is greater than $\, n \, $ and $\,
\theta >{1\over n}\, $ then $\, T_{\cal M}^{\theta }\,$ is reflexive.
Moreover, if the Cantor-Bendixson index of $\, \cal M \, $ is greater than $\,
\omega \, $ then $\, T_{\cal M}^{\theta }\,$ does not contain any $\, \ell ^p
\, ,$ while if the Cantor-Bendixson index of $\, {\cal M }\,$ is finite then
 $\, T_{\cal M}^{\theta }\,$ contains some $\, \ell ^p \, $ or $\, c_o \, .$ In
particular, if $\, {\cal M }=\{\,  A\subset {\bf N}\, :\, |A|\leq n\, \}\, $ and
$\, {1\over n}<\theta <1\,$ then $\, T_{\cal M}^{\theta }\,$ is isomorphic to
some $\, \ell ^p \, .$ } }

\medskip

{\bf Notation  }. We denote by $c_{00}$ the space of finitely supported
sequences.
  
$\{ e_n\} _{n=1}^{\infty}\,$ is the canonical basis in $c_{00}$.

 If $\, E\subset {\bf N}\,$ and $\,x\in c_{00}\,$, we denote by $\,Ex\,$ the
restriction of $x$ to $E$, that is $\, Ex=x. X _E$ .

 For $E\, ,F\, $ nonvoid subsets of ${\bf N}\,$ we write $\, E<F\, $ for $\,\max
 E<\min F\,$. We write $\, n<E\,$ for $\,\{ n\} <E$ .

If $\,E_i \subset {\bf N} ,\, i=1,\dots ,n\,$ and $\, E_1<E_2<\dots <E_n\,$ we
say that $\, E_1,\dots ,E_n \,$ are {\bf successive sets}.

For $\, x\in c_{00}\; ,\,  {\rm supp}(x)\,$ is the set $\, \{ i\, :\,x_i \neq
0\,\}$ .

If $\, x_i \in c_{00}\, ,\, x_i \neq 0\, ,\, i=1,\dots ,n\,$ and $\,
{\rm supp }(x_1)<\dots <{\rm supp}(x_n)\, $ we say that $\, x_1 ,\dots ,x_n \,$
are {\bf successive vectors}.

Finally, $\,\| \; .\; \| _p\, $ denotes the $\ell ^p$ norm, $\, 1\leq p\leq
\infty$.
\bigskip
We recall that the norm $\,\| \; .\; \| _T \,$ on Tsirelson's space $\,T\,$ is
defined by the following implicit equation.
$${\rm {For\, all}}\; x\in T\quad \| x \| _T =\max \{\, \| x\| _{\infty}\, ,\,
{1\over 2}\, \sup \sum _{i=1}^n \| E_i x\| _T\,\}\, ,$$
where the "sup" is taken over all $n$ and all sequences $\, ( E_1 ,\dots ,E_n
)\,$ of successive subsets of ${\bf N}$ with $\, n\leq \min E_1$.

It is well known that Tsirelson's space is a reflexive space not containing
any $\,\ell ^p\, $ ([T]).
\medskip
A basic notion in the definition of Tsirelson's space is that of a "Schreier
set". A finite set $A$ is called a Schreier set if $\, |A|\leq\min A .$ The
above mentioned properties of Tsirelson's space are consequences of the
following two properties of the class $\, \cal S\,$ of Schreier sets:

1)$\, \cal S\,$ does not contain infinite increasing sequences of sets.

2)$\, \cal S\,$ contains arbitrarily large (finite) sets.
\medskip
This leads us to the following generalization:

\medskip
Let $\, \cal M\,$ be a set of finite subsets of $\,{\bf N}.$ We say that
 $\, \cal M\,$ is {\bf compact} if the set 

\noindent $\, {\cal X} ({\cal M })\, =\{
X_A\, :\, A\in \, {\cal M}\,  \}\,$ is a compact subset of the Cantor space
$\, 2^{\bf N }.$

\bigskip
{\bf Definition 1} .{\it { Let $\, \cal M\,$ be as above. A family
 $\, ( E_1 ,\dots ,E_n)\,$ of successive finite subsets of $\,{\bf N}\, $is
said to be {\bf $\, \cal M$-admissible} if there exists a set $\, A=\{m_1
,\dots ,m_n\}\,\in  \cal M\,$ such that $\, m_1\leq E_1 <m_2 \leq E_2<\dots
<m_n\leq E_n\, .$}}
\bigskip
{\bf Definition 2} .{\it { Let $\, \cal M\,$ be a compact family. Let $\,
0<\theta <1\, .$ We define the $\, (\cal M ,\theta)$-Tsirelson-type space
$\,T_{\cal M}^{\theta}\, $ as follows:

 $\,T_{\cal
M}^{\theta}\, $ is the completion of $\,c_{00}\,$ under the norm $\,\|\,
.\,\|\,$ satisfying the implicit equation
$${\rm {For\, all}}\; x\in c_{00}\quad \| x \|  =\max \{\, \| x\| _{\infty}\,
,\, \theta \,\sup \sum _{i=1}^n \| E_i x\| \,\}\, ,$$
where the {\rm "sup"} is taken over all $n$ and all $\,\cal M $-admissible
families 
 $\, ( E_1 ,\dots ,E_n)\, .$}}
 \medskip
{\bf Remark }. It is easy to see that $\,\{e_n\}_{n=1}^{\infty}\,$ is a
1-unconditional basis for  $\,T_{\cal M}^{\theta}\, .$ 

\bigskip
Let us now recall the notion of the Cantor-Bendixson index of a countable
compact space $\, W\, .$

Given a countable compact space $\, W\, ,$ one defines the sequence $\,
W^{(\lambda )}\, ,\; \lambda <\omega _1\,$ as follows:
$$W^{(0)} =W$$
$$W^{(\lambda +1)}=\,\{\,x\in W\,:\, x {\rm {\, is \, a\, limit \, point \, of
\; }} W^{(\lambda )}\}$$
$${\rm {and\; for\; a\; limit\; ordinal\; \lambda\;}} , \quad W^{(\lambda )}
=\cap_{\mu <\lambda } W^{(\mu )}$$

The Cantor-Bendixson index of $\, W\,$ is then defined as the least $\,
\lambda < \omega _1\, $ for which $\, W^{(\lambda )}=\emptyset\, .$
\bigskip
We state without proofs the following results:
\medskip
{\bf Theorem 1} .{\it { Let $\,\cal M \,  , \theta\, $ be as in the definition
of
 $\,T_{\cal M}^{\theta}\, .$ Suppose that there is $\, n\in {\bf N}\, $ such
that the Cantor-Bendixson index of $\,\cal X (\cal M )\, $ is at least $\,
n+1\,$ and $\, {1\over n }<\theta <1\,.$ Then $\,T_{\cal M}^{\theta}\, $
is reflexive.}}
\medskip
{\bf Theorem 2} .{\it { If the Cantor-Bendixson index of $\,\cal X (\cal M )\, $
is greater than  $\,\omega \,$ then $\,T_{\cal M}^{\theta}\, $ does
not contain isomorphically any $\,\ell ^p\, .$}}
\bigskip
In the opposite direction we have the following result:
\medskip
{\bf Theorem 3} .{\it { If the Cantor-Bendixson index of $\, \cal X (\cal M
)\, $ is finite, then  the space $\,T_{\cal
M}^{\theta}\, $ contains an isomorphic copy of $\, \ell ^p \,,$
for some $\, 1<p<+\infty\, ,$ or $\, c_0 \, .$}}
\medskip
We shall present here a special case of Theorem 3, namely the case where

\noindent $\, {\cal M } =\{\, A\subset {\bf N}\, :\, |A|\leq n\,\}\, $ for some
$\, n\in {\bf N}\, .$
\medskip
{\bf Theorem 3a }.{\it { Let $\, n\in{\bf N}\, , n\geq 2 $ be fixed.
 If $\, {\cal M } =\{ \, A\subset {\bf N}\, :\, |A|\leq n\,\}\, $ and $\,
{1\over n}<\theta <1\,$ then the space $\,T_{\cal M}^{\theta}\, $ is isomorphic
to $\, \ell ^p\, ,$ where $\, {1\over p} +\log _n ({1\over \theta })=1\, .$

In other words:
The norm $\, \|\, ,\, \| \,$ defined on $\, c_{00}\, $ by the implicit equation

$$ \| x \|  =\max \{\, \| x\| _{\infty}\,
,\, {1\over {n^{1\over q} } }\, \sup \sum _{i=1}^d \| E_i x\| \,\}\, ,$$
where the "{\rm sup}" is taken over all $\, d\leq n\, $ and all sequences of
successive sets $\, (E_1 ,\dots E_d )\, ,$ is equivalent to the $\, \ell _p \,$
norm, where $\, {1\over p}+{1\over q} =1\, .$}}

\noindent {\bf Proof }: For the proof we first need to deal a little with the
dual of $\,\tmu\, .$

We define inductively a sequence $\, \{ K_s \} _{s=0} ^{\infty}\,$ of subsets
of $\, [-1\, ,\, 1]^{(N)} \, $ as follows:
$$K_0 =\{ {\buildrel + \over - }e_n \, :\, n\in {\bf N} \}$$

and for $\, s\in {\bf N}$
$$K_{s+1} =K_s \cup \{ \,\theta \cdot (f_1 +\dots +f_d \, )\, : \, d\leq n ,\,
f_i\, ,\i=1 \; {\rm {are\, successive\, and\, }}f_i \in K_s \, \} $$

It is not difficult to see that for every $\, x\in c_{00} \,$

$$\| x\| =\sup _{f\in \cup _{s=1} ^{\infty } K_s } <f\, ,\, x>$$

The proof of the Theorem goes through four steps:
\medskip
\noindent {\bf Step 1}

{\centerline {For every $\, x\in c_{00}\quad \| x \| \leq \| x\| _p .$ }} 

\noindent {\bf {Proof}} : It
is enough to show that for every $\, s\, $ and every
$\, f\in K_s \, , \, |f(x)|\leq \| x \| _p \, .$

This is done by induction on $\, s\, .$

For $\, f\in K_0 \,$ it is trivial. 
 Suppose that for some $\, s\, $ we have that for all $\, f\in K_s \,$ and all
$\, y\in c_{00}\quad |f(y)|\leq \| y \| _p \, .$
Let $\, x\in c_{00}\, $ and $\, f=\ncl (f_1 +\dots +f_d )\, \in K_{s+1}\, ,$
where $\, d\leq n\, , \,  f_1,\dots , f_n \, $ are successive and belong to
$\, K_s \,.$  

Then, setting $\, x_i = ({\rm supp}(f_i))(x) \, $ we have
$$|f(x)|\leq \ncl \sum _{i=1}^d |f_i (x)| =\ncl \sum _{i=1}^d |f_i (x_i)|\leq 
$$

$$\ncl \sum _{i=i}^d \| x_i\| _p\leq \big( {d\over n} \big) ^q
 (\sum _{i=i}^d \| x_i\| _p )^{{1\over p}} \leq \| x\| _p $$
by the induction hypothesis and H\"older's inequality.

\noindent {\bf Step 2} 

{\centerline {For all $\, m\in {\bf N} \quad {1\over {n^{1\over
p}} } m^{{1\over p}}\leq \nore .$}}
\medskip
\noindent {\bf {Proof}} : Suppose first
that $\, m=n^s \,$ for some $\, s\in {\bf
N}.\, $ The functional $\, f={1\over {n^{s\over q}}} \big( \sum_{i=1}^{n^s} e_i
\big)\,$ clearly belongs to $\, K_s \, .$ So
$$\| \sum_{i=1} ^{n^s} e_i \| \geq f\big( \sum_{i=1} ^{n^s} e_i \big) =
{1\over {n^{s\over q}}}n^s= n^{s(1-{1\over q} )}= n^{s\over
p}=m^{1\over p}$$

Now let $\, m\in {\bf N}\, $ and find $\, s\, $ such that $\, n^s \leq m\leq
n^{s+1} .\,$ Then 
$$\| \sum_{i=1} ^m e_i \| \geq \| \sum_{i=1} ^{n^s} e_i \| =n^{s\over
p}={1\over {n^{1\over
p}} }n^{{{s+1}\over p}}\geq {1\over {n^{1\over
p}} }m^{1\over p} .$$

\noindent {\bf Step 3}

For every normalized block sequence $\, (x_k)_{k=1} ^{\infty }\, $ of the
basis $\, (e_n )_{n=1} ^{\infty } \,$ we have
$$\| \sum a_k x_k \| \leq {2\over \theta } \| \sum a_k e_k \|$$

for all coefficients $\, (a_k)\, .$

\noindent {\bf {Proof }}: It is enough to show that for every $\, \phi \in \cup
_{s=1} ^{\infty } K_s \, $ one gets
$$\phi (\sum a_k x_k )\leq {2\over \theta } \| \sum a_k e_k \|$$
For the proof we need the following technical notions:
\medskip 
{\bf {Definition A}}.{\it { Let $\, m\in {\bf N}\, $ and $\, \phi \in K_m
\setminus K_{m-1}\, .$ We call {\bf analysis } of $\, \phi \, $ any sequence $\,
\{ F^s(\phi ) \}_{s=0}^m \,$ of subsets of  $\, \cup _{s=1}
^{\infty } K_s \, $ such that:

1) For every $\, s\;  F^s(\phi ) \,$ consists of successive elements  of $\,
K_s \;$ and

\noindent  $\cup _{f \in  F^s(\phi )}  {\rm supp}(f) ={\rm supp}(\phi )\, .$

2) If $\, f\, $ belongs to $\,  F^{s+1}(\phi ) \, $ then either
 $\, f\in  F^s(\phi )\, $ or there is  a $\, d\leq n \, $ and successive $\,
f_1 ,\dots ,f_d \in  F^s(\phi )\, $ with $\, f=\theta (f_1 +\dots +f_d )\, .$

3) $ F^m(\phi )=\{ \phi \} .$}}

\medskip
{\bf {Remark}}. It is clear by the definition of the sets $\, K_s \, $ that each
$\, \phi \in \cup K_s\, $ has an analysis.

Also one can check that if $\, f_1\in  F^s(\phi )\, ,\, f_2 \in  F^{s+1}(\phi
)\,$ then either $\, {\rm supp}(f_1)\subseteq  {\rm supp }(f_2)\,$ or $\,
{\rm supp }(f_1)\cap {\rm supp }(f_2) =\emptyset .$

\medskip
So let $\, \phi \in K_m \setminus K_{m-1}\, .$
By the 1-unconditionallity of $\, (e_k)\, $ we may and will assume that there is
$\, \ell \in {\bf N}\, $ such that $\, {\rm supp }(\phi ) =\cup _{k=1 } ^{\ell }
{\rm supp }(x_k )\, $ and the $\, x_k $ s and $\, \phi \,$ have only
non-negative coordinates.

\medskip
{\bf {Definition B}}. {\it {Let $\, \phi \, , \, (x_k )_{k=1} ^{\ell} \, $ be as
above. Let $\, \{  F^s(\phi ) \} _{s=0} ^m \, $ be a fixed analysis of $\,
\phi \, .$ For $\, k=1, \dots , \ell \, $ we set

$$s_k =\cases {   \max \{ \, s\, :\, 0\leq s <m {\rm {\; and \; there 
\; are \; at
\; least \; two\; }}f_1 ,\, f_2 \in  F^s(\phi ) \cr
\phantom {\max \{ \, }{\rm such \; that\; } f_i
(x_k)>0,\; i=1,\, 2\, \} , \cr
\phantom {\max }  {\it {\; when\; this\; set\; is\; non-empty }}\cr
\quad\cr
0\, , 
{\it \quad when\; } {\rm supp }(x_k) \; {\it { is \; a\; singleton\,
.}}}$$

So for each $\, k=1,\dots , \ell \,$ there exists a family $\, f_1 ^k ,\dots
,f_{d_k} ^k\; , \; 1\leq d_k \leq n \, $ of successive functionals in 
$\,  F^{s_k}(\phi )\, $ such that for each $\, i=1, \dots , d_k\, ,\;  f_i ^k
(x_k ) > 0 \, $ and 

\noindent $\, {\rm supp }(x_k) \subseteq \cup _{i=1 } ^{d_k } supp(f_i
)\, .$ 
We define the {\bf initial part } $\, x'_k\, $ and the {\bf final part}
$\, x''_k\, $ of $\, x_k \, $ with respect to $\, \{  F^s(\phi )\} _{s=0} ^m\,
$ as follows: 
$$x'_k =x_k | {\rm supp }(f_1 ^k )$$
$$x''_k =x_k |\cup _{i=2 } ^ {d_k} {\rm supp }(f_i ^k )$$
($\, x''_k =0\, $ if $\, {\rm supp } (x_k)\, $ is a singleton).}}
\medskip

Our aim is to show that
$$\phi (\sum _{k=1} ^{\ell} x'_k )\leq {1\over \theta}\| \sum _{k=1} ^{\ell}
a_k e_k \| $$
and
$$\phi (\sum _{k=1} ^{\ell} x''_k )\leq {1\over \theta}\| \sum _{k=1} ^{\ell}
a_k e_k \| $$
Since the proofs of these inequalities are similar we shall only prove the
first. In particular we shall show by induction on $\, s\leq m \, $ that for
every $\, J\subseteq \{ 1,\dots , \ell \} \,$ and every $\, f\in F^s(\phi ) $
$$(*)\qquad |f(\sum _{k\in J } a_k x_k )|\leq {1\over \theta }\| \sum _{k\in J}
a_k e_k \| \, .$$

$(*)\,$ is clear for $\, s=0\,.$
 Suppose that we know $\, (*)\,$ for $\, s<m\,$; we shall prove it for $\, s+1\,
.$

Let $\, J\subseteq \{1,\dots ,\ell \}.\,$ Let $\, f\in F^{s+1}(\phi )\, ,\;
f=\theta (f_1 +\dots +f_d ) \,$ with $\, d\leq n\, ,\, (fi)_{i=1} ^d\, $
successive members of $\, F^s(\phi )\, .$ Consider the sets 
$$K=\{ \, k\in J \,:\; {\rm there \; exists\; } i\in \{ 1,\dots ,d-1 \}\; {\rm
such \;that } 
\; f_i (x'_k )>0  \;{\rm and \;} f_{i+1}(x'_k )>0\, \}$$  
$$ I= \{ \, i\, :\; 1\leq i\leq d {\rm \; and \; there \; exists\; }K\in J
{\rm \; such \; that\; }{\rm supp }(x'_k) \subset {\rm supp }(f_i )\, \}$$
\noindent {\bf Claim } .  $|K|+|I|\leq n \, .$

\noindent {\bf Proof of the claim }. Let $\, k\in K\, .$ Then there exists $\,
i<d \,$ such that $\, f_i (x'_k )>0 \,$ and $\, f_{i+1} (x'_k)>0\, .$ This
means that $\, s<s_k \, ,$ so there exists $\, f\in F^{s_k}(\phi )\, $ with
$\, {\rm supp}(f_{i+1}) \subseteq {\rm supp}(f)\, .$ Then $\, f(x'_k) >0\, .$
But if $\, f\in F^{s_k}(\phi )\, $ and $\, f(x'_k )>0\, $ then by the
definition of $\, x'_k\, ,\; \max {\rm supp}(x'_k ) =\max {\rm supp}(f)\, .$ So
$\, {\rm supp}(f_{i+1})\subset {\rm supp} (x'_k)\, .$ Thus $\, i+1 \notin I\, .$
Therefore we can define a one-to-one map $\, G:\, K\to \{1, \dots , d \}
\setminus I \, ;$ hence $\, |K|\leq d-|I|\leq n-|I|\, .$
\medskip
We proceed with the proof of the inductive step. For $\, i\in I\, $ set

\noindent 
$E_i =\{\, k\in J\, :\, {\rm supp}(x'_k )\subseteq {\rm supp}(f_i)\, \}\, .$
Ofcourse $\, E_i \cap K =\emptyset \,$ for every $\, i\in I\, .$
We have
$$f(\sum _{k\in J}a_k x'_k)=\theta \sum _{k\in J}f_i (a_k x'_k)=
\theta \bigl\lbrack \sum _{i\in I} f_i (\sum _{k\in E_i} a_kx'_k)+\sum
_{k\in K} (\sum _{k=1} ^d f_i ) (a_k x'_k)\bigr\rbrack\, .$$
By the inductive hypothesis and the fact that for each $\, k\in J\, $

\noindent $\sum
_{k\in K} (\sum _{k=1} ^d f_i ) (a_k x'_k)\leq {1\over\theta }\| a_k x_k\|
={1\over\theta }\| a_k e_k \|\, ,$ we get  $$f(\sum _{k\in J}a_k x'_k)\leq
\theta \Big( {1\over \theta } \sum _{i\in I}\|E_i (\sum _{k\in J} a_k e_k )\|
+{1\over \theta }\sum _{k\in K} \| a_k e_k\| \leq {1\over \theta }\| \sum
_{k\in J} a_k e_k \| \Big) $$ using the fact that $\, |K|+|I|\leq n\,$ and the
implicit equation satisfied by the norm.

The proof of the inductive step and thus the proof of Step 3 are complete.

\noindent {\bf Step 4}

For all $\, \ell \, $ and all rational non-negative $\, (r_j )_{j=1} ^{\ell
}\,$

$$\| \sum _{j=1} ^{\ell } r_j ^{1\over p }e_j \| \geq {1\over {2n}} (\sum
_{j=1 }^{\ell } r_j )^{{1\over p}}$$
{\bf Proof} . Write $\, r_j ={{k_j}\over k}\, ,\, k_j, \, k \in {\bf N}\, .$
Set $\, s_0 =0\, ,\, s_j =k_1 +\dots +k_j\,  $ and 
$\,u_j =\sum _{i=s_{j-1}+1} ^{s_j} e_i \, ,$

\noindent $ j=1,\dots ,\ell\, .$ By Step 1,
$\, \| u_j\| \leq k_j ^{1\over p}\, .$ So
$$\| \sum _{j=1} ^{\ell } r_j ^{1\over p }e_j \| ={1\over {k^{{1\over
p}}}} \| \sum _{j=1} ^{\ell } k_j ^{1\over p }e_j \|\geq 
{1\over {k^{{1\over p}}}} \| \sum _{j=1}^{\ell} \| u_j \| e_j \|$$
by unconditionallity.

By Step 3
$${1\over {k^{{1\over p}}}} \| \sum _{j=1}^{\ell} \| u_j \| e_j \| \geq
{\theta \over 2}  {1\over {k^{{1\over p}}}} \| \sum _{j=1}^{\ell} \| u_j \|
{{u_j} \over {\| u_j \| }}\|=
{\theta \over 2}  {1\over {k^{{1\over p}}}}\| \sum _{j=1}^{\ell} \sum
_{i=s_{j-1}+1} ^{s_j} e_i\| ={\theta \over 2}  {1\over {k^{{1\over p}}}}\|
\sum _{i=1} ^{s_{\ell }}e_i \| $$
By Step 2, $\, \| \sum _{i=1} ^{s_{\ell }}e_i \| \geq {1 \over {n ^{{1\over
p}}}} s_{\ell }^{{1\over p }}\, $ so setting $\, \theta ={1\over {n ^{{1\over
q }}}}\, $ we get 
$${\theta \over 2}  {1\over {k^{{1\over p}}}}\| \sum _{i=1} ^{s_{\ell }}e_i
 \| \geq {1\over {2n}} {{s_{\ell}}^{{1\over p}}\over {k ^{{1\over p}}}}=
{1\over {2n}}\big( {{\sum _{j=1} ^{\ell }k_j } \over k }\big) ^{{1\over
p}}={1\over {2n}} \big( \sum _{j=1} ^{\ell }r_j \big) ^{{1\over p}} \, .$$ 
\medskip
Step 4 and the unconditionallity of $\, (e_n )_{n \in {\bf N}}\, $ imply that
$\, \| \sum a_k e_k \|\geq {1\over{2n}}(\sum |a_k |^p ) ^{{1\over p}}\, $
for all
coefficients $\, (a_k )\, .$ This fact combined with Step 1 completes the
proof of the Theorem.
\medskip {\bf Remark }. The result of the Theorem can also be deduced by Steps
1, 2 and 3 using a well known Theorem of Zippin [Z].
\medskip
\centerline {{\bf REFERENCES }}
\noindent [T]\quad  B. S. Tsirelson, Not every Banach space contains $\, \ell
^p\, $ or $\, c_o \, ,${\it { Funct. Anal. App. }} 8

\noindent \phantom{[T]}\quad (1974), p. 138-141

\noindent [Z]\quad M. Zippin, On perfectly homogeneous bases in Banach spaces,
{\it { Israel J. of Math. }}4 

\noindent \phantom{[Z]}\quad (1966), p. 265-272

\end